\documentclass[12pt,reqno]{amsart}
\usepackage{mathtools}
\usepackage{amssymb}
\usepackage{mathrsfs}
\usepackage{amsmath}
\usepackage{mathptmx}
\usepackage{comment}
\usepackage{enumerate}
\usepackage{stackrel}
\usepackage[dvipsnames]{xcolor}
\usepackage{amsopn}
\usepackage{dsfont}
\usepackage{hyperref}
\newcommand{\hide}[1]{}
\usepackage{enumerate}
\numberwithin{equation}{section}

\newtheorem{thm}{Theorem}

\usepackage{float}
\usepackage{graphicx}
\usepackage{caption}
\usepackage{subcaption}
\usepackage{titlesec}
\usepackage{fancyhdr}
\usepackage{textcomp}
\usepackage{fancyhdr}
\newtheorem{theorem}{Theorem}[section]

\newtheorem{corollary}[theorem]{Corollary}
\newtheorem{lemma}[theorem]{Lemma}

\newtheorem{notation}[theorem]{Notation}

\newtheorem{proposition}[theorem]{Proposition}
\newtheorem{definition}[theorem]{Definition}
\theoremstyle{definition}
\newtheorem{remark}[theorem]{Remark}

\DeclareMathOperator{\Capac}{cap}

\newcommand{\abs}[1]{| #1 |}

\DeclareMathOperator{\N}{\mathbb{N}}
\DeclareMathOperator{\Z}{\mathbb{Z}}
\DeclareMathOperator{\C}{\mathbb{C}}

\DeclareMathOperator{\T}{\mathbb{T}}

\DeclareMathOperator{\D}{\mathbb{D}}

\DeclareMathOperator{\IM}{Im}
\DeclareMathOperator{\h}{h}

\titleformat{\section}
  {\normalfont\scshape}{\thesection}{1em}{\centering}
\titleformat{\subsection}[runin]
  {\bfseries}{\thesubsection}{1em}{}
\DeclareMathOperator{\Hol}{Hol}

\usepackage{soul}
\usepackage[dvipsnames]{xcolor}

\begin{document}
\title{On the Hardy Number of Koenigs Domains}

\author[M.D. Contreras]{Manuel D. Contreras}
\address{Manuel D. Contreras. Departamento de Matem\'atica Aplicada II and IMUS, Escuela T\'ecnica Superior de Ingenier\'ia, Universidad de Sevilla, Camino de los Descubrimientos, s/n 41092, Sevilla, Spain}
\email{contreras@us.es}  

\author[F. J. Cruz-Zamorano]{Francisco J. Cruz-Zamorano}
\address{Francisco J. Cruz-Zamorano. Departamento de Matem\'atica Aplicada II and IMUS, Escuela T\'ecnica Superior de Ingenier\'ia, Universidad de Sevilla, Camino de los Descubrimientos, s/n 41092, Sevilla, Spain}
\email{fcruz4@us.es}  

\author[M. Kourou]{Maria Kourou}
\address{Maria Kourou. Julius-Maximilians-Universit\"at W\"urzburg, Institut f\"ur Mathematik, Emil Fischer Stra{\ss}e 40, 97074, W\"urzburg, Germany}
\email{maria.kourou@uni-wuerzburg.de}   

\author[L. Rodr\'iguez-Piazza]{Luis Rodr\'iguez-Piazza}
\address{Luis Rodr\'iguez-Piazza. Departmento de An\'alisis Matem\'atico and IMUS, Facultad de Matem\'aticas, Universidad de Sevilla, Calle Tarfia, s/n 41012 Sevilla, Spain}
\email{piazza@us.es}  

\thanks{M. D. Contreras, F. J. Cruz-Zamorano and L. Rodr\'iguez-Piazza are partially supported by Ministerio de Innovaci\'on y Ciencia, Spain, project PID2022-136320NB-I00, and Junta de Andaluc\'ia, project P20\_00664. }
\thanks{F. J. Cruz-Zamorano is also partially supported by Ministerio de Universidades, Spain, through the action Ayuda del Programa de Formaci\'on de Profesorado Universitario, reference FPU21/00258.}
\thanks{M. Kourou is partially supported by the Alexander von Humboldt Foundation.}

\fancyhf{}
\renewcommand{\headrulewidth}{0pt}
\fancyhead[RO,LE]{\small \thepage}
\fancyhead[CE]{\footnotesize M.D. CONTRERAS, F.J. CRUZ-ZAMORANO, M. KOUROU, AND L. RODR\'IGUEZ-PIAZZA}
\fancyhead[CO]{\footnotesize  ON THE HARDY NUMBER OF KOENIGS DOMAINS} 

\fancyfoot[L,R,C]{}

\begin{abstract}
This work studies the Hardy number for the class of hyperbolic planar domains satisfying Abel's inclusion property, which are usually known as Koenigs domains. More explicitly, we prove that for all regular domains in the above class, the Hardy number is greater or equal than $1/2$, and this lower bound is sharp. In contrast to this result, we provide examples of general domains whose Hardy numbers are arbitrarily small. Additionally, we outline the connection of the aforementioned class of domains with the discrete dynamics of the unit disc and obtain results on the range of Hardy number of Koenigs maps, in the hyperbolic and parabolic case. 
\end{abstract}

\maketitle

\section{Introduction}
Through the \textit{Hardy spaces} on the unit disc, $H^p(\D)$, one can define the \textit{Hardy number} of a holomorphic function $f \in \Hol(\D,\C)$, given by
$$\h(f) := \sup\left(\{ p>0 : f \in H^p(\D) \} \cup \{0\}\right) \in [0,+\infty],$$
which somehow measures the ``growth'' of $f$. Motivated by this concept, Hansen \cite{Hansen1970} introduced a similar idea in order to examine the range of holomorphic functions taking values in a domain $\Omega$. The so-called \textit{Hardy number} of a domain $\Omega$ is defined as
$$\h(\Omega) := \inf \{\h(f):  f\in \Hol(\D,\Omega)  \}. $$
The study of Hardy numbers focuses on unbounded domains, as for bounded domains $\Omega $ it follows trivially that $\h(\Omega) = +\infty$. 

A classical problem has been to determine the range of the Hardy number of a hyperbolic planar domain $\Omega$ (i.e., those whose boundary contains at least two points) in terms of its geometry and boundary behavior. As a matter of fact, several works have contributed to estimates on the Hardy number for certain classes of hyperbolic domains. Hansen \cite{Hansen1970, Hansen1971} provided a characterization of the Hardy number for starlike and spirallike with respect to the origin domains. A few years later, Ess\'en \cite{Essen_1981} obtained estimates of the Hardy number of a general hyperbolic domain in terms of harmonic measures and the logarithmic capacity. Thereafter, Kim and Sugawa \cite{Kim_Sugawa2011} examined the Hardy number for unbounded $K$-quasidisks. 
Quite recently, the range of the Hardy number for comb domains was studied by Karafyllia \cite{Kar_2022}.

The stepping stone, however, for the current work is the paper by Poggi-Corradini \cite{PC_1997_2}, who examined the Hardy number of hyperbolic domains $\Omega$ satisfying Schr\"oder's inclusion property; i.e. $\lambda \Omega \subseteq \Omega$, for some $\lambda \in \D$. Inspired by the impact of Schr\"oder's and Abel's functional equations on the discrete dynamics of the unit disc $\D$, our main focus lies on the range of the Hardy number for hyperbolic domains $\Omega$ satisfying \textit{Abel's inclusion property}; namely $\Omega+1\subseteq \Omega$. Due to its relevance, these domains will be called \textit{Koenigs domains}. 

One can easily notice that all domains in the aforementioned class are unbounded. Let us recall that a domain $\Omega$ in $\C$ is called {\sl regular} if the logarithmic capacity of $\C\setminus \Omega$ is positive.   If the domain $\Omega$ is non-regular, the Hardy number of $\Omega$ vanishes. This follows, for instance, from \cite[Theorems 5.1.1 and 5.4.2, p. 209 and 211]{Nevanlinna1970}, where it is proved that every universal covering map $p \colon \D \to \Omega$ has non-tangential limits almost nowhere on the unit circle $\T$ if the logarithmic capacity of $\C \setminus \Omega$ is zero. Our main result provides a lower bound for the Hardy number of a regular Koenigs domain:

\begin{theorem}
\label{thm:main-hardy}
Let $\Omega \subseteq \C$ be a Koenigs domain. Then, $\h(\Omega) \geq 1/2$ if and only if $\Omega$ is a regular domain. In particular, for every Koenigs domain $\Omega$ it holds that $\h(\Omega)\in \{0\}\cup [1/2,+\infty]$. 
\end{theorem}

The conclusion of the latter result does not hold for general domains, as shown in Section \ref{sec:small}, where we construct domains whose Hardy number is arbitrarily small. In fact, for every $p \in (0,+\infty)$, we provide a domain whose Hardy number is exactly $p$. 

The proof of Theorem \ref{thm:main-hardy} strongly depends on potential theory, namely a sharp estimate of the logarithmic capacity of compact sets obtained as union of integer translations of a compact non-polar subset of $\C$. Recall that a compact set in the complex plane is {\sl polar} if its logarithmic capacity is zero.

\begin{theorem}
\label{thm:main-potential}
Let $E$ be a compact non-polar subset of $\overline{D(0,1/4)}$. Set $$K_n = \bigcup_{j=1}^n(E+j).$$
Then, 
$$\lim_{n \to +\infty}\dfrac{\Capac(K_n)}{\Capac([0,n])} = 1.$$
Moreover, $\log(\Capac(K_n)) = \log(n/4) + \mathcal{O}(1/\sqrt{n})$.
\end{theorem}

After a preliminary introduction on logarithmic capacity, Hardy number and harmonic measure in Section \ref{sec:preliminaries}, we prove Theorem \ref{thm:main-potential} in Section \ref{sec:proof1} and Theorem \ref{thm:main-hardy} in Section \ref{sec:proof2}.
 After this, the aforementioned family of examples appears on Section \ref{sec:small}.

The main motivation for the latter results is that Abel's inclusion property is directly connected to non-elliptic discrete dynamics of the unit disc, mainly through the use of Koenigs maps. To this extent, it is of interest to study any implications of Theorem \ref{thm:main-hardy}in the context of discrete iteration theory. To motivate this, several examples are presented in Section \ref{sec:app} to relate the Hardy number of Koenigs maps and the properties of its associated self-map. We conclude by discussing some implications of this work for inner functions in Section \ref{sec:inner}. 

\begin{notation}
We write $a_n = \mathcal{O}(b_n)$ if there exists $C > 0$ and $N \in \N$ such that $\abs{a_n} \leq Cb_n$ for all $n \geq N$.
\end{notation}

\section{Preliminaries}
\label{sec:preliminaries}
\subsection{Logarithmic capacity.}
\label{subsec:log}

One of the main tools in potential theory  is the potential of a measure, which is defined as follows:
\begin{definition}
\label{def:potential-energy}
Let $\mu$ be a finite positive measure on $\C$ with compact support. Its (logarithmic) potential is the function $p_{\mu} \colon \C \to [-\infty,+\infty)$ given by
$$p_{\mu}(z) = \int_{\C}\log\abs{z-w}d\mu(w), \quad z \in \C.$$
Its (logarithmic) energy $I(\mu) \in [-\infty,+\infty)$ is given by
$$I(\mu) = \int_{\C}p_{\mu}(z)d\mu(z) = \int_{\C}\int_{\C}\log\abs{z-w}d\mu(w)d\mu(z).$$
From these concepts one can define the (logarithmic) \textit{capacity} of a set $X \subseteq \C$, that is
\begin{equation}\label{eq:capacitydef}
    \Capac(X) = \sup_{\mu}\exp(I(\mu))
\end{equation}
where the supremum is taken among every probability measure $\mu$ whose support is a compact subset of $X$.
\end{definition}

In these definitions, we understand that $I(\mu) = -\infty$ if the former integral is not convergent. Indeed, if $\Capac(X) = 0$, then $I(\mu) = -\infty$ for every probability measure $\mu$ with compact support in $X$. If this happens, $X$ is said to be a \textit{polar set}.
A property is said to be satisfied \textit{nearly everywhere} in $X \subset \mathbb{C}$, if it is satisfied for all points in $X$, except maybe for a Borel polar subset. In general, polar sets are negligible from the potential-theoretic point of view. Moreover, the following lemma holds:
\begin{lemma}
\cite[Corollary 3.2.5]{Ransford}
\label{lem:countable}
A countable union of Borel polar sets is polar.
\end{lemma}

Let $\Omega \subseteq \C$ be a domain, that is, an open and connected set. If $\Capac(\C \setminus \Omega) > 0$, then $\Omega$ is called \textit{regular}. Furthermore, when it comes to compacta, the theory of the potentials is richer: if $X\subset \C$ is a compact set, then there exists a probability measure $\nu$ with support in $X$ attaining the supremum in \eqref{eq:capacitydef}. In fact, this measure is unique if $X$ is non-polar and in such a case $\nu$ is said to be the \textit{equilibrium measure} of $X$; see \cite[Theorem 3.7.6]{Ransford}. The potential associated to the equilibrium measure satisfies the following property:
\begin{thm}[Frostman's Theorem]
\cite[Theorem 3.3.4]{Ransford}
\label{thm:Frostman}
Let $X \subseteq \C$ be a non-polar compact set, and let $\nu$ be its equilibrium measure. Then, $p_{\nu}(z) \geq I(\nu)$ for all $z \in \C$. Moreover, $p_{\nu} \equiv I(\nu)$ nearly everywhere on $X$.
\end{thm}

Further details on regular domains, logarithmic capacity and polar sets may be found in \cite[Chapters 3-5]{Ransford}.

\subsection{Hardy numbers.}
\label{subsec:hardy}
Given $0 < p < +\infty$, the Hardy space $H^p(\D)$ is defined as the set of all holomorphic maps $f \colon \D \to \C$ such that
$$\sup_{0<r<1}\int_{\T}\abs{f(r\xi)}^pdm(\xi) < +\infty,$$
where $m$ is the normalized length measure in the boundary of the unit disc $\T$. In the case $p = +\infty$, the Hardy space $H^{\infty}(\D)$ stands for the set of all holomorphic maps $f \colon \D \to \C$ that are bounded, that is,
$$\sup_{z \in \D}\abs{f(z)} < +\infty.$$
We refer to \cite{Duren_Hp} for a complete introduction to this topic. One of the properties of these spaces is that they form a decreasing family, that is, $H^p(\D) \supseteq H^q(\D)$ if $0 < p \leq q \leq +\infty$. These relations suggest the following idea: given a holomorphic map $f \colon \D \to \C$, its Hardy number $\h(f)$ is defined as
$$\h(f) = \sup(\{0\} \cup \{p > 0 : f \in H^p(\D)\}) \in [0,+\infty].$$
Note that $f \in H^p(\D)$ for every $0 < p < \h(f)$, and $f \not\in H^p(\D)$ if $p > \h(f)$.

In a similar manner, this idea can be translated to domains. Given a domain $\Omega \subseteq \C$, its Hardy number is defined as
$$\h(\Omega) = \inf\{\h(f) : f \in \text{Hol}(\D,\Omega)\}.$$
Here, $\text{Hol}(\D,\Omega)$ denotes the set of all holomorphic maps $f \colon \D \to \C$ such that $f(\D) \subseteq \Omega$. Let us state some well-known properties of the Hardy number of a domain:
\begin{lemma}
\cite[Lemmas 2.1 and 2.3]{Kim_Sugawa2011}
\label{lem:properties}
Let $\Omega, \Omega^{\prime} \subseteq \C$ be two domains. Then 
\begin{enumerate}[\normalfont(a)]
\item $\h(\Omega)=+\infty$, if $\Omega$ is bounded.
\item $\h(\Omega^{\prime}) \leq \h(\Omega)$, if $\Omega\subseteq \Omega^{\prime}$.
\item $\h(\varphi(\Omega))=\h(\Omega)$ for a complex affine map $\varphi(z)=az+b$, $a\neq 0$. 
\item $\h(\Omega)=0$, if $\mathbb{C} \setminus \Omega$ is bounded.
\item $\h(\Omega) \geq \frac{1}{2}$ if $\Omega \neq \C$ is simply connected.
\item $\h(\Omega) = \h(p)$, where $p \colon \D \to \Omega$ is a universal covering map of $\Omega$.
\end{enumerate}
\end{lemma}

The Hardy number for certain simply connected domains is already known. For instance, $\h(H)=1$ for a half-plane $H$. In the case of a strip e.g. $S(a,b)=\{ z \in \mathbb{C} :  a< \IM z <b\}$, with $a,b \in \mathbb{R}$, it is known that $\h(S(a,b))= +\infty$. The Hardy number of sectors has also been examined in \cite{Hansen1970}: suppose $\theta \in (0,2\pi]$ and 
$S_{\theta}: = \left\{re^{i \phi} :r>0, |\phi| <\frac{\theta}{2}\right\}$, then $\h(S_{\theta}) =\frac{\pi}{\theta}$.

A characterization of the Hardy number for certain domains can be obtained through the following classical result:
\begin{thm}
\label{thm:Nevanlinna}
\cite[Theorems 5.1.1 and 5.4.2, p. 209 and 211]{Nevanlinna1970}
A universal covering map $p \colon \D \to \Omega$ has non-tangential limits almost everywhere on $\T$ if and only if $\Omega$ is a regular domain.
\end{thm}
Recall that any function whose Hardy number is positive (that is, it belongs to some Hardy space) has non-tangential limits almost everywhere on $\T$. Then, from the above result and the properties of the Hardy numbers, it follows that the Hardy number of every non-regular domain is zero.

\subsection{Harmonic measure.} For a given regular domain $\Omega \subseteq\C$, let $B \subseteq\partial \Omega$ be a Borel set. The harmonic measure $\omega(z,B,\Omega)$ of $B$ at a point $z \in \Omega$ is the solution of the generalized Dirichlet problem in $\Omega$ with boundary values $1$ on $B$ and $0$ on $\partial \Omega \setminus B$. 

For a fixed Borel set $B \subseteq\partial \Omega$, $\Omega \ni z \mapsto \omega (z, B, \Omega)$ is a harmonic and bounded function. In addition, for a fixed point $z \in \Omega$, the map $B \mapsto \omega(z, B, \Omega)$ is a Borel probability measure on $\partial \Omega$. We refer to \cite[Section 4.3]{Ransford} for an introduction to harmonic measure.

In \cite[Lemma 1]{Essen_1981}, Ess\'en proposed a relation between Hardy number and harmonic measure, which was later improved by Kim and Sugawa in the following result
\begin{thm}
\label{thm:KimSugawa}
\cite[Lemma 3.2]{Kim_Sugawa2011}
Let $\Omega \subseteq\C$ be a domain with $0 \in \Omega$. Then,
$$\h(\Omega) = \liminf_{R \to + \infty}\left(-\dfrac{\log \omega(0,F_R,\Omega_R)}{\log R}\right),$$
where $\Omega_R$ is the connected component of $\Omega \cap D(0,R)$ containing the origin and $F_R = \partial\Omega_R \cap \{|z| = R\}$.
\end{thm}
This result will be useful in Section \ref{sec:proof2} to deduce Theorem \ref{thm:main-hardy} from Theorem \ref{thm:main-potential}, and in Section \ref{sec:small} to provide examples of domains $\Omega$ satisfying $\h(\Omega) = p$ for every given $p \in (0,+\infty)$.

\section{Proof of Theorem \ref{thm:main-potential}}\label{sec:proof1}

Before moving on to the proof of Theorem \ref{thm:main-potential}, we prove two auxiliary lemmas related to the construction of the equilibrium measure of the compact sets $K_n$, $n \in \mathbb{N}$.

In the sequel we will use the equilibrium measure of a compact interval, which can easily be derived from \cite[Eq. (1.7), p. 25]{SaffTotik}. For any $n \in \N$, the equilibrium measure $\mu$ of the interval $[0,n]$ is absolutely continuous (with respect to Lebesgue's measure $m$) and it is given by
\begin{equation}
\label{eq:mu}
\dfrac{d\mu}{dm}(t) = \dfrac{\chi_{[0,n]}(t)}{\pi \sqrt{t(n-t)}}.
\end{equation}
In particular, it follows (cf. \cite[Eq. (1.8), p. 25]{SaffTotik})
\begin{equation}
\label{eq:properties_mu}
p_{\mu}(t) = I(\mu) = \log(\Capac([0,n])) = \log(n/4), \quad \text{for all } t \in [0,n]\, ; 
\end{equation}
the above result is a combination of Theorem \ref{thm:Frostman} and \cite[Theorems 4.2.2 and 4.2.4]{Ransford}.

Using this equilibrium measure, we define the following coefficients:
$$\alpha_{j}:= \int_{j-1}^jd\mu(t) = \frac{1}{\pi} \int_{j-1}^{j} \frac{dt}{\sqrt{t(n-t)}}, \quad j \in \N, \quad 1 \leq j \leq n.$$
Notice that $\alpha_j$ also depends on $n$. However, in the seek of clearance, this it not explicitly written in the notation. Some properties of these numbers can be derived directly from the definition. For example, it follows that $\alpha_j>0$ for all $j=1,..., n$ and that $\sum_{j=1}^n \alpha_{j}=1$. It is also possible to notice that these numbers are endowed with some symmetry, namely $\alpha_j = \alpha_{n-j+1}$, and that $\alpha_j$ is non-increasing for $j=1, ..., \lfloor (n+1)/2 \rfloor$, where $\lfloor x \rfloor$ denotes the integer part of the real number $x$.

We prove the following estimations:

\begin{lemma}\label{lem:properties_aj}
\begin{enumerate}[\normalfont(a)]
\item There exists $C_1 > 0$ such that $\alpha_j \leq C_1/\sqrt{n}$, for all $n \in \N$ and all $j = 1,\ldots,n$.
\item There exists $C_2 > 0$ such that $$\sum_{\substack{j=1 \\ j\neq k}}^n \frac{\alpha_j}{|j-k|} \leq \frac{C_2}{\sqrt{n}},$$ for all $n \in \N$ and all $k = 1,\ldots,n$.

\end{enumerate}
\begin{proof}
(a) Fix a natural number $n\geq 2$. Then 
$$\alpha_1 = \dfrac{1}{\pi} \int_0^1\dfrac{dt}{\sqrt{t(n-t)}} \leq \dfrac{1}{\pi\sqrt{n-1}}\int_0^1\dfrac{dt}{\sqrt{t}} = \mathcal{O}\left(\dfrac{1}{\sqrt{n}}\right).$$
Due to the monotonicity and symmetric properties of the coefficients $\alpha_j$, we have that $\alpha_j\leq \alpha_1$, for $j=1, ...,n$, so that (a) holds.

(b) Fix $n \in \N$ and define
$$S(k) = \sum_{\substack{j=1 \\ j\neq k}}^n \frac{\alpha_j}{|j-k|}, \quad k = 1,\ldots,n.$$
By symmetry, notice that $S(k) = S(n-k+1)$. Thus, it is enough to work with $k \in \N$ such that $k \leq (n+1)/2$. If this is the case, notice that for any $j \in \N$ with $1 \leq j \leq (n+1)/2$ it is possible to check that $\alpha_j = \alpha_{n-j+1}$ but $\abs{j-k} \leq \abs{n-j+1-k}$. Therefore,
$$S(k) \leq 2 \sum_{\substack{j=1 \\ j\neq k}}^{\lfloor (n+1)/2 \rfloor} \frac{\alpha_j}{|j-k|}.$$

Given two finite sequences $\{ a_1, a_2, ..., a_m \}$ and $\{ b_1, b_2, ..., b_m \}$ of non-negative real numbers, the Hardy-Littlewood inequality asserts that
$$
\sum _{j=1}^m a_j b_j\leq \sum _{j=1}^m a_j^* b_j^*,
$$
where $\{a_j^*\}$ denote the sequence of elements $a_j$ arranged in decreasing order and similarly for  $\{b_j^*\}$ (see \cite[\S 10.2]{inequalities} and also \cite[p. 43]{Bennett-Colin}).
Notice that the map $\{1,\ldots,\lfloor (n+1)/2 \rfloor\} \ni j \mapsto \alpha_j$ is decreasing. Moreover, taking $b_j = 1/|j-k|$, for $j=1, ...,  \lfloor (n+1)/2 \rfloor$, $j\neq k$, and $b_{k}=0$, we have that $b_j^*\leq 2/j$. Thus
$$\sum_{\substack{j=1 \\ j\neq k}}^{\lfloor (n+1)/2 \rfloor} \frac{\alpha_j}{|j-k|} =
 \sum_{j=1}^{\lfloor (n+1)/2 \rfloor} \alpha_j b_{j} \leq 2\sum_{j=1}^{\lfloor (n+1)/2 \rfloor}\dfrac{\alpha_j}{j}.$$

But now, notice that, for $n\geq 2$,
\begin{equation*}
\begin{split} \sum_{j=1}^{\lfloor (n+1)/2 \rfloor}\dfrac{\alpha_j}{j}&  \leq 
\dfrac{2}{\pi}\int_0^{(n+1)/2}\dfrac{dt}{(t+1)\sqrt{t(n-t)}} \leq \dfrac{2}{\pi\sqrt{(n-1)/2}}\int_0^{+\infty}\dfrac{dt}{(t+1)\sqrt{t}} = \mathcal{O}\left(\dfrac{1}{\sqrt{n}}\right), 
\end{split}
\end{equation*}
where, in the first inequality, we have used that $1/j \leq 2/(t+1)$ whenever $j-1\leq t\leq j$. Thus, the result follows.
\end{proof}
\end{lemma}

To proceed with the proof of Theorem \ref{thm:main-potential}, let $E \subseteq \overline{D(0,1/4)}$ be a compact non-polar set with $0\in E$ and equilibrium measure $\nu$. Fix $n \in \mathbb{N}$ and set
$$K_n=\bigcup_{j=1}^n (E+j).$$
Let us define the positive measure $\sigma$ given by
$$\sigma(A) = \sum_{j=1}^{n} \alpha_j \nu(A-j),$$
where $A \subseteq \C$ is a Borel set. Notice that $\sigma$ is a probability measure whose support is a compact set lying on $K_n$.
Once more, observe that $\sigma$ depends on $n$.

Recalling the Definition \ref{def:potential-energy}, let us prove the following lemma:
\begin{lemma}\label{lemma:sigma}
Under the above notation, the following statements hold:
\begin{enumerate}[\normalfont(a)]
\item $p_{\sigma}(x) \geq \log(n/4) + \mathcal{O}(1/\sqrt{n})$ for every $x \in K_{n}$.
\item $p_{\sigma}(x) = \log(n/4) + \mathcal{O}(1/\sqrt{n})$ for nearly every $x \in K_n$.
\item $\abs{p_{\sigma}(x)-p_{\sigma}(y)} = \mathcal{O}(1/\sqrt{n})$ for nearly every $x,y \in K_n$.
\item $p_{\sigma}(x) = \log(n/4) + \mathcal{O}(1/\sqrt{n})$ for every $x \in E$.
\item $\abs{p_{\sigma}(x)-p_{\sigma}(y)} = \mathcal{O}(1/\sqrt{n})$ for every $x \in E$ and nearly every $y \in K_n$.
\end{enumerate} 
In fact, the underlying constants do not depend on $x$ and $y$.
\begin{proof} In this proof we will always assume that $n\geq 2$. Note that, by the definition of the probability measure $\sigma$, its potential $p_{\sigma}$ can be written as
$$p_{\sigma}(z) = \sum_{j=1}^n  \alpha_j p_{\nu}(z-j) \, , \quad z \in \C.$$

Let $k \in \{1,...,n\}$. Notice that, by Theorem \ref{thm:Frostman},  $p_{\nu}(x-k) \geq  I(\nu)$ for all $x\in E+k$ and $p_{\nu}(x-k) =  I(\nu)$ nearly everywhere in $x\in E+k$. Therefore \begin{align*}
p_{\sigma}(x) & = \sum_{j=1}^n \alpha_j p_{\nu}(x-j) \geq \alpha_k I(\nu) + \sum_{\substack{j=1 \\ j\neq k}}^n \alpha_j \int_E \log|x-j-y| d\nu(y)
\end{align*}
for all $x\in E+k$ and the inequality turns to be an equality  nearly everywhere. We claim that, given $x\in E+k$, 
\begin{equation} \label{eq:claim}
\alpha_k I(\nu) + \sum_{\substack{j=1 \\ j\neq k}}^n \alpha_j \int_E \log|x-j-y| d\nu(y)= \log\left(\dfrac{n}{4}\right) + \mathcal{O}\left(\dfrac{1}{\sqrt{n}}\right).
\end{equation}
Thus, assuming this claim,  we clearly conclude (a) and (b).

Let us prove the claim. Take $x\in E+k$ and write $x^{\prime} =x-k \in E$. Then
\begin{align*}
\log|x-j-y| & = \log\abs{x^{\prime} -y + k-j} = \log\abs{k-j} + \log\left| 1+ \frac{x^{\prime} -y}{k-j}\right| \\
& = \log\abs{k-j} + \mathcal{O}\left(\frac{1}{\abs{k-j}} \right),
\end{align*}

for all $y \in E$, since $\abs{x^{\prime}-y} \leq 1/2$. 
Due to Lemma \ref{lem:properties_aj}.(b), we obtain 
\begin{equation}\label{eq:psigma}
\alpha_k I(\nu) + \sum_{\substack{j=1 \\ j\neq k}}^n \alpha_j \int_E \log|x-j-y| d\nu(y) = \alpha_k I(\nu) + \sum_{\substack{j=1 \\ j \neq k}}^n\alpha_j  \log|k-j| + \mathcal{O} \left(\frac{1}{\sqrt{n}} \right) \,  .
\end{equation}

Using the equilibrium measure $\mu$ for $[0,n]$ as described in \eqref{eq:mu}, let us define
\begin{align*}
S & := \alpha_k I(\nu) + \sum_{\substack{j=1 \\ j \neq k}}^n  \alpha_j\log|k-j|  - p_{\mu}(k-1/2) \\
& =  \alpha_k I(\nu) + \sum_{\substack{j=1 \\ j \neq k}}^n \alpha_j\log|k-j|  - \sum_{j=1}^n \int_{j-1}^j \log \left| k-1/2-t \right| d\mu(t) \\
& = \alpha_k I(\nu) - \int_{k-1}^k \log \left|k-1/2-t \right|d\mu(t) + \sum_{\substack{j=1 \\ j \neq k}}^n  \int_{j-1}^j  \left(\log\abs{k-j} - \log \left|k-1/2-t \right|\right)d\mu(t).
\end{align*}
We want to prove  that $\abs{S} = \mathcal{O}(1/\sqrt{n})$. This can be done by examining each of the three terms that define $S$. By Lemma \ref{lem:properties_aj}.(a), since $\abs{I(\nu)} < +\infty$, one concludes that $\abs{\alpha_k I(\nu)} = \mathcal{O}(1/\sqrt{n})$.

Concerning the second term, let us define
$$F(k) = \int_{k-1}^k \log \left|k-1/2-t \right|d\mu(t), \quad k = 1, \ldots, n.$$
Using the symmetry, notice that $F(n-k+1) = F(k)$. Therefore, we can suppose that $1 \leq k \leq (n+1)/2$. On the one hand, if $k = 1$, then
\begin{align*}
\left|\frac{1}{\pi}\int_{0}^1 \frac{\log \abs{1/2 -t}}{\sqrt{t(n-t)}} dt\right| & \leq
\frac{1}{\pi}\int_{0}^1 \frac{\big|\log \abs{1/2 -t}\big|}{\sqrt{t(n-t)}} dt \\
& \leq \frac{1}{\pi\sqrt{n-1}} \int_{0}^1 \frac{\big|\log \abs{1/2 -t}\big|}{\sqrt{t}}dt =
\mathcal{O}\left(\dfrac{1}{\sqrt{n}}\right).
\end{align*}

On the other hand, if $2 \leq k \leq (n+1)/2$, then
\begin{align*}
\left|\frac{1}{\pi} \int_{k-1}^k\frac{\log \abs{k-1/2-t}}{\sqrt{t(n-t)}}dt \right|& \leq
\frac{1}{\pi} \int_{k-1}^k\frac{\big|\log \abs{k-1/2-t}\big|}{\sqrt{t(n-t)}}dt \\
& \leq \frac{1}{\pi\sqrt{(k-1)(n-k)}}\int_{k-1}^k \big|\log \abs{k-1/2-t} \big| \, dt \\
& \leq \frac{1}{\pi\sqrt{(n-1)/2}}\int_0^1 \big|\log \abs{t-1/2} \big| \, dt
= \mathcal{O}\left(\dfrac{1}{\sqrt{n}}\right).
\end{align*}

Let us now estimate the third term. Fix $j = 1, \ldots, n$ with $j \neq k$, and notice that all $t \in [j-1,j]$ can be written as $t = j-1+t'$ for some $t' \in [0,1]$. Therefore,
$$\log\abs{k-j} - \log \abs{k-1/2-j+1-t'} = \log \left|1+\dfrac{1-2t'}{2(k-j)}\right| = \mathcal{O}\left(\dfrac{1}{\abs{k-j}}\right),$$
where we have used that $\left|\frac{1-2t'}{2(k-j)}\right|\leq \frac12$.
Thus, there exists $c > 0$ such that
\begin{align*}
\sum_{\substack{j=1 \\ j \neq k}}^n  \int_{j-1}^j  \left(\log\abs{k-j} - \log \left|k-1/2 -t \right|\right)d\mu(t) & \leq \sum_{\substack{j=1 \\ j \neq k}}^n  \int_{j-1}^j  \dfrac{c}{\abs{k-j}}d\mu(t) \\
& = \sum_{\substack{j=1 \\ j \neq k}}^n  c\dfrac{\alpha_j}{\abs{k-j}} = \mathcal{O}\left(\dfrac{1}{\sqrt{n}}\right),
\end{align*}
where Lemma \ref{lem:properties_aj}.(b) has been used.

Combining these arguments, we conclude that $\abs{S} = \mathcal{O}(1/\sqrt{n})$. As a result, from the definition of $S$ and \eqref{eq:psigma}, we obtain that
\begin{equation}
\label{eq:potentialK}
\alpha_k I(\nu) + \sum_{\substack{j=1 \\ j\neq k}}^n \alpha_j \int_E \log|x-j-y| d\nu(y)= p_{\mu}(k-1/2) + S + \mathcal{O}\left(\dfrac{1}{\sqrt{n}}\right) = \log\left(\dfrac{n}{4}\right) + \mathcal{O}\left(\dfrac{1}{\sqrt{n}}\right).
\end{equation}
where \eqref{eq:properties_mu} has been used, what proves the claim and thus (a) and (b).

It is clear that (c) follows from (b).

Furthermore, (d) follows from similar ideas to the ones used to obtain (b). To see this, notice that $\abs{x-y} \leq 1/2$ for all $x,y \in E$, from which it follows that, for any given $k \in \N$,
\begin{equation}
\label{eq:log}
\log\abs{x-y-k}-\log(k) = \log \left| 1-\dfrac{x-y}{k}\right| = \mathcal{O}\left(\dfrac{1}{k}\right).
\end{equation}
Therefore, for every $x \in E$,
\begin{equation}
\label{eq:psigmaE}
p_{\sigma}(x) = \sum_{k = 1}^n\alpha_k\int_E\log\abs{x-y-k}d\nu(y) = \sum_{k = 1}^n\alpha_k\log(k) + \mathcal{O}\left(\dfrac{1}{\sqrt{n}}\right),
\end{equation}
where \eqref{eq:log} and Lemma \ref{lem:properties_aj}.(b) have been used.

From \eqref{eq:properties_mu}, $p_{\mu}(0)=\log(n/4)$ and we have 
\begin{align*}
S & := \sum_{k = 1}^n\alpha_k\log(k) - \log(n/4) = \sum_{k = 1}^n\alpha_k\log(k) - p_{\mu}(0) \\
& = \sum_{k = 1}^n\int_{k-1}^k\log(k)d\mu(t) - \sum_{k = 1}^n\int_{k-1}^k\log(t)d\mu(t)  = \sum_{k = 1}^n\int_{k-1}^k\log(k/t)d\mu(t).
\end{align*}

We claim that $\abs{S} = \mathcal{O}(1/\sqrt{n})$. To see this, notice that if $k = 1$, then
$$\int_{0}^1\log(1/t)d\mu(t) = \dfrac{1}{\pi}\int_0^1\dfrac{-\log(t)}{\sqrt{t(n-t)}}dt \leq \dfrac{1}{\pi\sqrt{n-1}}\int_0^1\dfrac{-\log(t)}{\sqrt{t}} dt = \mathcal{O}\left(\dfrac{1}{\sqrt{n}}\right).$$
On the other hand, if $k \geq 2$, notice that
$$\log\left(\dfrac{k}{k-1}\right) = \log\left(1+\dfrac{1}{k-1}\right) = \mathcal{O}\left(\dfrac{1}{k}\right).$$
Therefore, there exists some $c > 0$ such that the following holds:
\begin{align*}
\sum_{k=2}^n\int_{k-1}^k\log(k/t)d\mu(t) & \leq \sum_{k=2}^n\int_{k-1}^k\log(k/(k-1))d\mu(t) \\
& \leq \sum_{k=2}^n\int_{k-1}^k\dfrac{c}{k}d\mu(t) = c\sum_{k=2}^n\dfrac{\alpha_k}{k} = \mathcal{O}\left(\dfrac{1}{\sqrt{n}}\right),
\end{align*}
by Lemma \ref{lem:properties_aj}.(b). Notice that the claim follows from these arguments. As a result, from the definition of $S$ and \eqref{eq:psigmaE}, we obtain that
$$p_{\sigma}(x) = \log\left(\dfrac{n}{4}\right) + \mathcal{O}\left(\dfrac{1}{\sqrt{n}}\right).$$

Lastly, (e) follows from (b) and (d).
\end{proof}
\end{lemma}

At this point, we can proceed with the proof of Theorem \ref{thm:main-potential}.

\begin{proof}[\bf Proof of Theorem \ref{thm:main-potential}]
Consider the equilibrium measure $\beta$ for $K_n$. From Theorem \ref{thm:Frostman}, $I(\beta)=p_\beta(x)$ for nearly every $x \in K_n$. Moreover, notice that $\sigma (A)=0$ for every Borel subset $A$ with $\Capac(A)=0$, as $\sigma$ depends on a finite translation of the equilibrium measure of $E$. Then,
$$\log(\Capac(K_n)) = \int_{K_n}p_{\beta}(x)d\sigma(x) = \int_{K_n}p_{\sigma}(x)d\beta(x),$$
where Fubini's Theorem has been used. Then, by Lemma \ref{lemma:sigma}.(c),
$$\abs{\log(\Capac(K_n)) - p_{\sigma}(y)} = \left|\int_{K_n}(p_{\sigma}(x)-p_{\sigma}(y))d\beta(x)\right| \leq \int_{K_n}\abs{p_{\sigma}(x)-p_{\sigma}(y)}d\beta(x) = \mathcal{O}\left(\dfrac{1}{\sqrt{n}}\right)$$
for nearly every $y \in K_n$. Thus, the result follows from Lemma \ref{lemma:sigma}.(b) and \eqref{eq:properties_mu}.
\end{proof}

\section{Proof of Theorem \ref{thm:main-hardy}} \label{sec:proof2}
\begin{proof}[\bf Proof of Theorem \ref{thm:main-hardy}]Using Theorem \ref{thm:Nevanlinna}, we will restrict to the case where $\Omega$ is a regular Koenigs domain.

Consider $X = \C \setminus \Omega$, and define $B_N = \{z \in X : \abs{z} \leq N\}$ for all $N \in \N$. Notice that $\cup_{N \in \N}B_N = X$. Then, by Lemma \ref{lem:countable}, since $\Capac(X) > 0$, there must exists some $N \in \N$ such that $\Capac(B_N) > 0$. Moreover, given $z \in B_N$, define $E_z = \{w \in B_N : \abs{w-z} \leq 1/4\}$. Notice that $\cup_{z \in B_N}E_z = B_N$, which is a compact set. Then, there must exist a (finite) sequence $ \{x_n \}_{n=1,...,m}$ in $ B_N$, such that $\cup_{n =1}^{m}E_{x_n} = B_N$. In that case, using Lemma \ref{lem:countable} once more and the fact that $\Capac(B_N) > 0$, there must exist $x \in B_N$ such that $E_x \subseteq X$ with $\Capac(E_x) > 0$.

Without loss of generality, we can suppose that $x = 0$. By construction, $E_x \subseteq \overline{D(0,1/4)}$. With this in mind, define $E = -E_x$, and consider
$$\Omega^* = \C \setminus \bigcup_{n \in \N} (E+n).$$
Note that $\Omega \subseteq -\Omega^*$, which yields $\h(\Omega) \geq \h(-\Omega^*) = \h(\Omega^*)$. Therefore, it is enough to prove that $\h(\Omega^*) \geq 1/2$. To do so, Theorem \ref{thm:main-potential} is used to estimate $\h(\Omega^*)$ in the following way:

As in Theorem \ref{thm:KimSugawa}, consider the harmonic measure $\omega_R(z) = \omega(z,F^*_R,\Omega^*_R)$, where $\Omega^*_R$ is the connected component of $\Omega^* \cap D(0,R)$ containing the origin and $F^*_R = \partial\Omega^*_R \cap \{w \in \C : |w| = R\}$. For a given $R > 3$, define $m = \max\{j \in \N : j \leq R/3\}$, and notice that
$$\Omega^*_R \cup \{w \in \C : m + 1/4 < \abs{w} < R\}) = D(0,R) \setminus K_{m}=: \Delta_R,$$
where $K_{m}:= \bigcup_{j = 1}^m(E+j)$. Therefore, consider $\beta_R(z) = \omega(z,F_R^*,\Delta_R)$ and $\tilde \beta_R(z) = \omega(z,\partial D(0,R),\Delta_R)$. By the maximum principle for harmonic functions \cite[Theorem 1.1.8]{Ransford} it is possible to show that $\omega_R(0) \leq \beta_R(0) \leq \tilde \beta_R(0)$. 

Furthermore, $\tilde \beta_R$ can be estimated through the use of potentials. To do so, we follow the ideas which were developed in Section \ref{sec:proof1}. Consider $\nu$ the equilibrium measure for $E$, and define the probability measure given by
$$\sigma(A) = \sum_{j = 1}^m\alpha_j\nu(A-j),
$$
for every Borel set $A \subset \C$. Notice that if $z \in \C$ is such that $\abs{z} = R$, then $p_{\sigma}(z) \geq \log(R-m-1/4) \geq \log(R/3)$ if $R$ is big enough. Take $\lambda_{m}=\inf\{p_{\sigma}(x):\, x\in K_{m}\}$. As $\sigma$ is supported by $K_{m}$, by the Minimum Principle for logarithmic potentials (see \cite[Theorem 3.1.4]{Ransford}), $\lambda_{m}=\inf\{p_{\sigma}(x):\, x\in \C\}$. Using both statements (a) and (b) of  Lemma \ref{lemma:sigma}, we have that 
$$\lambda_{m}=\log(m/4)+\mathcal{O}\left(\dfrac{1}{\sqrt{m}}\right) =p_{\sigma}(y)+ \mathcal{O}\left(\dfrac{1}{\sqrt{m}}\right) ,
$$ 
for nearly every $y \in K_m$. Therefore, the function $\gamma_{R}:=p_{\sigma}-\lambda_{m}$ is non-negative on $\C$, harmonic in $\C\setminus K_{m}$, and for $z \in \C$ such that $\abs{z} = R$ we have
\begin{equation*}
\begin{split}
\gamma_{R}(z)&\geq \log(R/3)-\lambda_{m}=\log(R/3)-\log(m/4)+\mathcal{O}\left(\dfrac{1}{\sqrt{m}}\right)\\
&\geq \log(R/3)-\log(R/12)+\mathcal{O}\left(\dfrac{1}{\sqrt{m}}\right)=\log(4)+\mathcal{O}\left(\dfrac{1}{\sqrt{m}}\right) = \log(4)+\mathcal{O}\left(\dfrac{1}{\sqrt{R}}\right).
\end{split}
\end{equation*}
Therefore, $\gamma_{R}(z)\geq 1$ whenever $|z|=R$ for $R$ large enough.
Thus $\tilde \beta_R(0)\leq \gamma_{R}(0)=p_{\sigma}(0)-\lambda_{m}$, for $m$ (and then $R$) large enough.

By Lemma \ref{lemma:sigma}.(e),  for nearly every $y \in K_m$,
$$
\gamma_{R}(0)=p_{\sigma}(0)-\lambda_{m}=p_{\sigma}(0)-p_{\sigma}(y)+ \mathcal{O}\left(\dfrac{1}{\sqrt{m}}\right) =\mathcal{O}\left(\dfrac{1}{\sqrt{m}}\right) =\mathcal{O}\left(\dfrac{1}{\sqrt{R}}\right).
$$

With this in mind, Theorem \ref{thm:KimSugawa} yields
$$\h(\Omega^*) = \liminf_{R \to + \infty}\left(-\dfrac{\log\omega_R(0)}{\log R}\right) \geq \liminf_{R \to + \infty}\left(-\dfrac{\log\gamma_R(0)}{\log R}\right) \geq \dfrac{1}{2}.$$
\end{proof}
\medskip

\section{Domains with arbitrary prefixed Hardy number}
\label{sec:small}
In order to strengthen the scope of Theorem \ref{thm:main-hardy} it is convenient to provide examples of general domains whose Hardy numbers are arbitrarily small. However, as far as we know, the vast literature on this topic is mainly devoted to simply connected domains for which the Hardy number is always greater or equal than $1/2$. We are not aware of any previous example whose Hardy number is in $(0,1/2)$. Therefore, to complete the work, we show that for every $p \in (0,1/2)$ there exists a domain $\Omega \subseteq \C$ such that $\h(\Omega) = p$. This construction involves a continuity argument, Theorem \ref{thm:KimSugawa}, and the maximum principle for harmonic functions \cite[Theorem 1.1.8]{Ransford}.

Before the construction, we start by proving the following lemma, where the notation $C_r = \{z \in \C : \abs{z} = r\}$ for $r > 0$ is used:
\begin{lemma}
\label{lemma:compact}
Let $K \subseteq \C$ be a non-polar compact set such that $0$ is in the unbounded connected component of $\C \setminus K$. Take $R_0 > 0$ with $K \subseteq \overline{D(0,R_0)}$. The function
$$R \mapsto \omega_R := \omega(0,C_R,D(0,R) \setminus K)$$
is continuous in $(R_0,+\infty)$.
\begin{proof}
First of all, using the maximum principle for harmonic functions, note that $\omega_R$ is a decreasing function of $R$, and so it is enough to prove that
$$\lim_{r \to r_0^-}\omega_r \leq \omega_{r_0} \leq \lim_{r \to r_0^+}\omega_r, \quad r_0 \in (R_0,+\infty).$$

For the first inequality, let $R_0 < r \leq r_0$, and define $u_r(z) = \omega(z,C_r,D(0,r) \setminus K)$. In particular, $\omega_r = u_r(0)$. Let us pick $r_1 \in (R_0,r_0)$ and note that $u_{r_0}$ is uniformly continuous in $\{z \in \C : r_1 \leq \abs{z} \leq r_0\}$. In particular, for every $\epsilon > 0$ there exists $\rho \in (r_1,r_0)$ such that $1-\epsilon \leq u_{r_0}(z) \leq 1$ if $\rho \leq \abs{z} \leq r_0$. By the maximum principle, if $r_1 \leq r \leq r_0$, then $(1-\epsilon)u_r \leq u_{r_0}$ in $\{z \in \C : r_1 \leq \abs{z} \leq r_0\}$. In this case, $\lim_{r\to r_0^-}(1-\epsilon)\omega_r \leq \omega_{r_0}$. Taking $\epsilon \to 0^+$, we are done.

Similarly, for the second inequality, let $R_0 < r_0 \leq r$ and define
$$v_r(z) = \omega(z,C_r,\{R_0 < \abs{w} < r\}) = \dfrac{\log(\abs{z}/R_0)}{\log(r/R_0)};$$
see \cite[Table 4.1]{Ransford}. By the maximum principle, if $R_0 < \abs{z} < r$, then $u_r(z) \geq v_r(z) = v_r(\abs{z})$, which means that $u_r(z)/v_r(r_0) \geq 1$ if $\abs{z} = r_0$. In particular, by the maximum principle, $u_{r_0} \leq u_r/v_r(r_0)$ in $D(0,r_0) \setminus K$. In this case, $\omega_{r_0} \leq \lim_{r \to r_0^+}\omega_r/v_r(r_0)$, where $\lim_{r \to r_0^+}v_r(r_0) = 1$, and so the lemma follows.
\end{proof}
\end{lemma}
To present the examples, we introduce the following notation: given $p > 0$ and a strictly increasing sequence $\{R_n\}$ in $(1,+\infty)$ such that $R_n \to + \infty$ as $n \to + \infty$, define the domain
\begin{equation}
\label{eq:small}
\Omega = \C \setminus \left(\bigcup_{n = 1}^{\infty}(C_{R_n} \setminus \Gamma_n)\right), \quad \Gamma_n = \left\lbrace z = R_ne^{i\theta} : \abs{\theta} < \dfrac{\pi}{R_n^{2p}}\right\rbrace.
\end{equation}
Let us also define
$$\Omega_{R,N} = D(0,R) \setminus \left(\bigcup_{n = 1}^N (C_{R_n} \setminus \Gamma_n) \right), \quad \omega_{R,N} = \omega\left(0, C_R, \Omega_{R,N}\right), \quad R > R_N.$$
We can now prove the following:
\begin{theorem}
For any $p \in (0,+\infty)$ there exists a domain $\Omega \subseteq \C$ for which $\h(\Omega) = p$.
\begin{proof}
Fix $p > 0$. The domain $\Omega$ is constructed as in \eqref{eq:small}, where the sequence $\{R_n\}$ is inductively constructed depending on $p$. To do this, set $R_1 = 2$. Suppose that $R_n$ is constructed for all $1 \leq n \leq N$. Let $R > R_N$, and notice that the maximum principle implies that
\begin{equation}
\label{eq:limomega}
\omega_{R,N} \leq 
\omega(0,\Gamma_N,D(0,R_N)) = \dfrac{1}{R_N^{2p}},
\end{equation}
as $\omega(0,\cdot,D(0,R))$ is the normalized Lebesgue measure on $C_R$. In particular, $\lim_{R \to R_N^+}\omega_{R,N} \leq 1/R_N^{2p}$, where the limit converges because the maximum principle implies that the function $R \mapsto \omega_{R,N}$ is decreasing for $R > R_N$.

We now claim that there must exists $R > R_N$ such that $\omega_{R,N} = 1/R^p$. To see this, arguing by contradiction, using the continuity given in Lemma \ref{lemma:compact} and the above inequality, we have that $\omega_{R,N} < 1/R^p$ for all $R > R_N$. Then, Theorem \ref{thm:KimSugawa} implies that $\h(\C \setminus (\cup_{n=1}^N(C_{R_n} \setminus \Gamma_n)) \geq p$. But this contradicts Lemma \ref{lem:properties}.(d). Thus, it is possible to define $R_{N+1} = \min\{R > R_N : \omega_{R,N} = 1/R^p\}$. From \eqref{eq:limomega}, we see that $R_{N+1} \geq R_N^2$, which yields that $R_N \to + \infty$ as $N \to + \infty$.

It is clear that for every $R > R_1$ there exists $N \in \N$ with $R_N < R \leq R_{N+1}$. From the construction, it follows
$$\omega(0,C_R, \Omega \cap D(0,R)) = \omega_{R,N} \leq 1/R^p,$$
where equality is attained for $R = R_{N+1}$. Therefore, applying Theorem \ref{thm:KimSugawa}, $\h(\Omega) = p$.
\end{proof}
\end{theorem}
\begin{remark}
In the setting of Lemma \ref{lemma:compact}, it is possible to prove that $\lim_{R \to + \infty} \omega_R\log(R) > 0$. This fact could be used instead of Theorem \ref{thm:KimSugawa} to show the claim we made in order to prove the existence of $R_{N+1}$.
\end{remark}

\section{Applications to Discrete Iteration Theory}
\label{sec:app}
Theorem \ref{thm:main-hardy} can be of interest in the non-elliptic discrete dynamics of the unit disc. To introduce this, let $\phi$ be a holomorphic self-map of $\D$. Let us denote by $(\phi_n)_{n\in \mathbb{N}}$, where ${\phi_n= \phi \circ ... \circ \phi}$ composing $n$ times, the \textit{sequence of iterates} of $\phi$.

As usual, if $\phi$ has a fixed point in $\D$, then $\phi$ is called \textit{elliptic}. In the case where $\phi$ is not an elliptic automorphism of $\D$, the Denjoy--Wolff Theorem asserts that $(\phi_n)$ converges locally uniformly in $\D$ to a point $\tau_{\phi} \in \overline{\D}$. The point $\tau_{\phi}$ is called the \textit{Denjoy--Wolff point} of $\phi$.

The research spectrum of the current work focuses on the case where $\phi$ has no fixed points in $\D$ and thus, $\tau_{\phi} \in \partial \D$. For this case, it is known that $0 < \phi^{\prime}(\tau_{\phi}) \leq 1$ (in the angular limit sense). More specifically, if $\phi^{\prime}(\tau_{\phi}) <1$, then $\phi$ is called \textit{hyperbolic}, whereas, if $\phi^{\prime} (\tau_{\phi})=1$, $\phi$ is called \textit{parabolic}.

Parabolic self-maps are further divided into two subcategories depending on its \textit{hyperbolic step}. 
The parabolic self-map $\phi:\D\to \D$ is \textit{of zero hyperbolic step} if for some --and hence, for all-- $z \in \D$, it is true that $\rho_{\D}(\phi_{n}(z), \phi_{n+1}(z)) \to 0$, $n \to +\infty$, where $\rho_{\D}$ denotes the pseudo-hyperbolic distance in $\D$. Otherwise, $\phi$ is parabolic \textit{of positive hyperbolic step}. For a further introduction on this topic, we refer to a recent book by Abate \cite[Chapter 4]{AbateBook}.

A prominent tool in examining properties of the iterates $(\phi_n)$ are the solutions to the so-called Abel's equation, that is, holomorphic maps $\sigma \colon \D \to \C$ satisfying
$$\sigma \circ \phi = \sigma + 1.$$
In the course of the past century, the existence of solutions for this equation in the case where $\phi$ is non-elliptic was shown. Indeed, this was achieved through the individual works by several mathematicians: Valiron examined the hyperbolic case \cite{Valiron}, while the two different parabolic cases were covered by Pommerenke \cite{Pom_1979,BakerPommerenke}, the second one in collaboration with Baker. The main idea of these three works was to explicitly construct a solution $\sigma$ to the Abel's equation for a given self-map $\phi$ in terms of its ``normalized'' iterates.

Some time afterwards, Cowen \cite{Cowen} proposed to find solutions of a general functional equation related to the self-map $\phi$. This eventually led to the following concept: a triple $(\Omega_0, \Psi, \Phi)$ is said to be a \textit{model} for $\phi$ if $\Omega_0\subseteq \C$ is a domain, $\Psi \colon \D \to \Omega_0$ is a holomorphic function, and $\Phi$ is an automorphism of $\Omega_0$, satisfying the following three conditions:
\begin{enumerate}[(i)]
\item \begin{equation}\label{eq:Schr/Abel}
\Psi \circ \phi = \Phi \circ \Psi,
\end{equation}
\item $$\bigcup_{n \in \mathbb{N}} \Phi^{-n}(\sigma(\D)) = \Omega_0,$$
\item there exists a domain $A \subseteq \D$ such that $\Psi$ is injective on $A$, $\phi(A) \subseteq A$, and for every $z \in \D$ there exists $n \in \N$ with $\phi_n(z) \in A$. 
\end{enumerate}

The domain $\Omega_0$ is called a \textit{base space}, the holomorphic function $\Psi$ is called an \textit{intertwining map}, the automorphism $\Phi$ is called a \textit{normal form}, and $A$ is called an \textit{absorbing set}.

Cowen showed that every non-elliptic self-map $\phi$ admits a model, which might be chosen so that it follows a certain canonical form. This can be formulated in the following way, where the notation $\mathbb{H} = \{z \in \C : \IM(z) > 0\}$ is used:
\begin{thm}\label{thm:models}
Let $\phi$ be a non-elliptic holomorphic self-map of $\D$. There exists a model $(\Omega_0, \sigma, \Phi)$ for $\phi$, where $\Phi(z) = z+1$, so that:
\begin{enumerate}[(i)]
\item $\phi$ is hyperbolic if and only if there exists $ \lambda > 1$ such that $\Omega_0 = \mathbb{S}(\lambda) := \{z \in \C : 0 < \IM(z) < \pi/\log(\lambda)\}$.
\item $\phi$ is parabolic of positive hyperbolic step if and only if $\Omega_0 = \mathbb{H}$ or $\Omega_0 = -\mathbb{H}$.
\item $\phi$ is parabolic of zero hyperbolic step if and only if $\Omega_0 = \C$.
\end{enumerate}
\end{thm}
\begin{remark}
This theorem is stated in \cite[Theorem 4.6.8]{AbateBook}, although it has been reformulated in such a way that \eqref{eq:Schr/Abel} always corresponds to the Abel's equation for $\phi$.
\end{remark}

Every non-elliptic holomorphic self-map admits an essentially unique holomorphic model (see \cite[Corollary 3.5.9]{AbateBook}). As a byproduct, the function $\sigma$ in Theorem \ref{thm:models} is unique if we assume that $\text{Re}(\sigma(0)) = 0$ in cases (i) and (ii), and $\sigma(0) = 0$ in case (iii). From now on, such normalized solution to the Abel's equation $\sigma$ will be called the \textit{Koenigs map} of $\phi$. The image $\Omega = \sigma(\D)$ plays a prominent role in the properties of $\phi$. Note that, from Abel's equation, $\Omega$ needs to satisfy that $\Omega + 1 \subseteq \Omega$. Therefore, $\Omega$ is a \textit{Koenigs domain} in the sense of the previous sections. Indeed, the conclusion in Theorem \ref{thm:main-hardy} might be strengthened if $\Omega$ is assumed to be a regular Koenigs domain of a certain self-map $\phi$ of $\D$, depending on the properties of $\phi$.

For example, if $\phi$ is hyperbolic, its Koenigs domain $\Omega$ is contained in a horizontal strip. This means that this Koenigs domain and the associated Koenigs map have an infinite Hardy number, see Subsection \ref{subsec:hardy}. Therefore, in the hyperbolic case, the Koenigs map is in $H^p$ for all $0 < p < +\infty$. A little more can be obtained through the space BMOA of analytic functions with bounded mean oscillation, that is, the set of analytic functions $f \colon \D \to \C$ with
$$\sup_{w \in \D}\|f_w\|_2 < +\infty, \quad \text{where} \quad  f_w(z) = f\left(\dfrac{z+w}{1+\overline{w}z}\right)-f(w), \quad z \in \D.$$
We refer to \cite{Girela} for a complete introduction to BMOA. It holds $H^{\infty} \subseteq \text{BMOA} \subseteq H^p$, for $ 0 < p < +\infty$. The key point is that Riemann mappings from $\D$ onto any strip $\mathbb{S}(\lambda)$ are known to be in BMOA. Therefore, by a subordination argument (see \cite[Corollary 10.1]{Girela}), the Koenigs map $\sigma$ corresponding to a hyperbolic self-map $\phi$ is also in BMOA. This can also be noticed geometrically using a deeper result of Hayman and Pommerenke, see \cite[Theorem 1]{HaymanPommerenke}.

In the parabolic case, $\Omega$ characterizes the hyperbolic step of $\phi$, as noticed in Theorem \ref{thm:models}. Therefore, if $\phi$ is a parabolic self-map of positive hyperbolic step, the conclusion on Theorem \ref{thm:main-hardy} can be improved: its Koenigs domain $\Omega$ is contained in a half-plane, and so the Hardy number of the Koenigs domain and the associated Koenigs map is always greater that one.

Furthermore, for every $p \in [1,+\infty]$, it is possible to find a parabolic self-map $\phi$ of positive hyperbolic step such that $\h(\Omega) = \h(\sigma) = p$ for some Koenigs map $\sigma$. In the case where $p$ is finite, this construction is achieved by means of sectors: fix some $\theta \in (0,\pi]$ and set the domain $S(\theta):=\{z \in \mathbb{C}: \text{arg}(z) \in (0,\theta)\}$. Let $\sigma \colon \D \to S(\theta)$ be a Riemann mapping, and define a self-map $\phi \colon \D \to \D$ given by $\phi(z) = \sigma^{-1}(\sigma(z)+1)$. By Theorem \ref{thm:models}, $\phi$ is a parabolic self-map of positive hyperbolic step and, up to a translation, $\sigma$ is the Koenigs map for $\phi$. Therefore, the properties of Hardy numbers show that $\h(\sigma) = \h(\sigma(\D)) = \h(S(\theta)) = \pi/\theta$; see Subsection \ref{subsec:hardy}. By choosing $\theta \in (0,\pi]$ appropriately, it is clear that for any $p \in [1,+\infty)$ one can find explicit examples with $\h(\sigma) = p$.

Similarly, for the infinite value of Hardy number, we consider the domain $\Omega = \{z =x+iy \in \mathbb{C} : x> 0, \,  0<y< \sqrt{x} \}$. With the definitions given above, note that for every $\theta \in (0,\pi/2)$ there exists $t \geq 0$ such that $\Omega+t \subseteq S(\theta)$. Due to the properties of the Hardy numbers, it follows $\h(\Omega)=\h(\Omega+t) \geq \h(S(\theta)) = \pi/\theta$. Letting $\theta \to 0^+$ leads to $\h(\Omega)=+\infty$. Once more, consider a Riemann mapping $\sigma \colon \D \to \Omega$ and define a self-map $\phi \colon \D \to \D$ given by $\phi(z) = \sigma^{-1}(\sigma(z)+1)$. By Theorem \ref{thm:models}, $\phi$ is a parabolic self-map of $\D$ of positive hyperbolic step and, up to a translation, $\sigma$ is a Koenigs map for $\phi$. Once more, the properties of Hardy numbers lead to $\h(\sigma) = \h(\sigma(\D)) = \h(\Omega) = +\infty$.

While considering parabolic self-maps of zero hyperbolic step, Koenigs domains are not included in some horizontal half-plane of $\C$. Therefore the main difference is that now it may actually occur that $\C \setminus \Omega$ has zero logarithmic capacity and thus, $\h(\sigma)$ can be zero. For example, set $\Omega = \C \setminus \{n \in \Z : n < 0\}$ and let $p \colon \D \to \Omega$ be a universal covering map with $p(0) = 0$. Using \cite[Theorems 8.1 and 8.2]{CDMP2007}, one can see that there exists a parabolic inner function $\phi \colon \D \to \D$ of zero hyperbolic step for which $p$ is its Koenigs map in the sense of Theorem \ref{thm:models}, but $\h(p) = \h(\Omega) = 0$.

In the regular case, the latter situation can not happen since Theorem \ref{thm:main-hardy} can be applied. It can also be shown that no uniform upper bound for $\h(\sigma)$ nor $\h(\Omega)$ exists. To be more precise: for every $p \in [1/2,+\infty]$ there exists a parabolic self-map $\phi$ of zero hyperbolic step such that $\h(\sigma) = \h(\Omega) = p$ for the Koenigs map for $\phi$. Similarly as before, in the finite cases,  this is done by means of the sectors $S(\theta):=\left\{ z\in \C : |\text{arg}(z)| < \theta/2 \right\}$ with $\h(S(\theta)) = \pi/\theta$ and a suitable Riemann mapping. The argument is also similar for the infinite value, where $\Omega = \{z=x+iy \in \C : x > 0, \, |y| < \sqrt{x}\}$ with $\h(\Omega) = + \infty$ is used.

To sum up, the foregoing considerations can be abbreviated as follows:
\begin{proposition}
\label{prop:sumup}
Let $\phi$ be a holomorphic non-elliptic self-map of $\D$ and consider its associated Koenigs map $\sigma$.
\begin{itemize}
\item[(i)] If $\phi$ is hyperbolic, $\sigma$ is in BMOA.
\item[(ii)]  If $\phi$ is parabolic of positive hyperbolic step, then $\h(\sigma) \geq 1$. Indeed, for every $p \in [1,+\infty]$ there exists a parabolic self-map $\phi$ of positive hyperbolic step for which $\h(\sigma) = p$.
\item[(iii)] If $\phi$ is parabolic of zero hyperbolic step, then $\h(\sigma) \in \{0\} \cup [1/2,+\infty]$. Indeed, for every ${p\in[1/2,+\infty]}$ there exists a parabolic self-map $\phi$ of zero hyperbolic step for which $\h(\sigma) = p$. Moreover, there exists an inner parabolic self-map of zero hyperbolic step, for which $\h(\sigma) = 0$.
\end{itemize}
\end{proposition}

\section{Dynamics of Inner Functions}\label{sec:inner}
Inspired by Proposition \ref{prop:sumup}.(iii), we end this paper with a brief result on the dynamics of inner functions. According to a classical result of Fatou, a self-map $\phi$ of $\D$ possesses radial limits almost everywhere on $\T$. Hence the \textit{boundary map} $\phi \colon \T \to \C$ defined by
$$\phi(\zeta) :=\lim_{r\to 1} \phi(r \zeta)$$
is well-defined for almost every $\zeta \in \T$. 

A classical, yet challenging, question in the discrete dynamical setting concerns boundary convergence for the sequence of iterates of inner parabolic self-maps of zero hyperbolic step, which is also studied in \cite[Chapter 4]{DoerMane1991}. Boundary convergence is strongly related to the boundary behavior of the associated Koenigs map, as noticed by the first named author, D\'iaz-Madrigal and Pommerenke. To show this, they used the construction of the Koenigs map that was developed in \cite{Pom_1979,BakerPommerenke}. However, the existence of an absorbing set can be used to relate the latter construction with the one given in Theorem \ref{thm:models}; see \cite[Theorem 3.1]{CDMPom_abel}. Taking this into account, the following is true:

\begin{thm}\cite[Theorem 1.1]{CDMPom_2012}\label{thm:boundconv}
Let $\phi$ be a parabolic self-map of $\D$ of zero hyperbolic step with Denjoy--Wolff point $\tau_{\phi}$. If the associated Koenigs map $\sigma$ has an angular limit almost everywhere on $\T$, then $(\phi_n(\xi))$ converges to $\tau_{\phi}$ for almost every $\xi \in \T$.   
\end{thm}

The latter result motivates the following:
\begin{theorem}\label{thm:inner}
Let $\phi:\D \to \D$ be a parabolic self-map with associated Koenigs map $\sigma$. 
\begin{enumerate}[\normalfont(i)]
\item If $\sigma$ has angular limits almost nowhere on $\partial \D$, then $\phi$ is inner.
\item If $\phi$ is inner and of zero hyperbolic step, then $\sigma$ has an angular limit either almost everywhere or almost nowhere on $\partial \D$.
\end{enumerate}
\begin{proof}
(i) We follow a similar technique as in the proof of \cite[Theorem 8.1]{CDMP2007}. Set $A:=\{\zeta \in \T: \angle \lim_{z\to \zeta} \sigma(z) \in \mathbb{C} \}$ and $B:=\left\{\xi \in \T \, : \, \angle\lim_{z \to \zeta} \phi(z)\in \D \right\}$. Since the Koenigs map $\sigma$ satisfies Abel's functional equation, it follows that for all $\xi \in B$,
$$\angle \lim_{z\to \xi} \sigma(z) = \angle \lim_{z\to \xi} (\sigma(\phi(z)) -1) = \sigma (\phi(\xi)) -1 \in \mathbb{C}.$$
As a result $B \subseteq A$. Therefore, if $A$ has zero measure, so it does $B$, meaning that $\phi$ is an inner function.  

(ii) With the notation from above, let us now suppose that $\phi$ is an inner function. In that case, for all $\xi \in A$ it follows that
$$\angle \lim_{z \to \xi } \sigma(\phi(z)) =\angle \lim_{z \to \xi } \sigma(z) +1 \in \mathbb{C}$$
and so, $\phi(\xi) \in A$, which implies $\phi(A) \subseteq A$. Conversely, if $\phi(\xi)\in A$ for some $\xi \in \T$, then 
$$\angle \lim_{z \to \xi } \sigma(z)  =\angle \lim_{z \to \xi }  \sigma(\phi(z)) -1 \in \mathbb{C},$$
which implies $\phi^{-1} (A) := \{\xi \in \T : \phi(\xi) \in A\} \subseteq A$. As a result, $A=\phi^{-1}(A)$.

Moreover, for an inner self-map $\phi$ of $\D$, its corresponding boundary map is ergodic if and only if $\phi$ has zero hyperbolic step (see e.g. \cite[Theorem 3.1]{DoerMane1991}, \cite[\S\ 5.1]{BMS_2005} or \cite[Chapter 6]{Aaronson_1997}). Therefore, from $A = \phi^{-1}(A)$, ergodicity implies that $A$ has either zero or full measure, meaning that $\sigma$ has an angular limit either almost everywhere or almost nowhere on $\T$. 
\end{proof}
\end{theorem}

As a matter of fact, concerning inner self-maps of $\D$ and taking Theorem \ref{thm:inner} into account, we prove that the latter condition on the boundary behavior of the Koenigs map can be weakened.
\begin{corollary}\label{cor:Boundconv_final}
Let $\phi$ be an inner parabolic self-map of $\D$ of zero hyperbolic step with Denjoy--Wolff point $\tau_{\phi}$ and associated Koenigs map $\sigma$. If $\sigma$ has non-tangential limit on a set of positive (Lebesgue) measure, then $\phi_n \to \tau_{\phi}$ a.e. on $\partial \D$.
\begin{proof}
As proved in Theorem \ref{thm:inner}, if $\sigma$ has angular limits on a set of positive measure, then it has angular limits almost everywhere on $\T$. Therefore, the result follows from Theorem \ref{thm:boundconv}.
\end{proof}
\end{corollary}

\end{document}